\documentclass[11pt]{amsart}
\input{diagrams}

\usepackage{bibentry}
\usepackage{amsmath}
\usepackage{amsthm}
\usepackage{amssymb}
\usepackage{amsfonts}
\usepackage{amsxtra}     
\usepackage{amscd}
\usepackage{epsfig}
\usepackage{verbatim}
\usepackage{latexsym,amstext,epsfig}
\usepackage[all, knot]{xy}
\xyoption{arc}

\newsymbol\pp 1275

\newcommand{\Hom}{\operatorname{Hom}}
\newcommand{\End}{\operatorname{End}}

\newcommand{\Rep}{\operatorname{Rep}}
\newcommand{\SI}{\operatorname{SI}}
\newcommand{\SL}{\operatorname{SL}}
\newcommand{\GL}{\operatorname{GL}}
\newcommand{\ZZ}{\mathbb Z}
\newcommand{\CC}{\mathbb C}
\newcommand{\RR}{\mathbb R}

\newcommand{\Ima}{\operatorname{Im}}

\newcommand{\tr}{\operatorname{Tr}}

\newtheorem{theorem}{Theorem}[section]
\newtheorem{proposition}[theorem]{Proposition}
\newtheorem{corollary}[theorem]{Corollary}
\newtheorem{lemma}[theorem]{Lemma}

\theoremstyle{definition}

\newtheorem{remark}[theorem]{Remark}

\newtheorem{example}[theorem]{Example}

\title[]{Quivers, long exact sequences and Horn type inequalities II}

\author{Calin Chindris}
\address{School of Mathematics, University of Minnesota, Minneapolis, MN, USA}
\email{chindris@math.umn.edu}

\begin{document}
\bibliographystyle{plain}
\begin{abstract}
We study the set of all $m$-tuples $(\lambda(1),\dots,\lambda(m))$
of possible types of finite abelian $p$-groups $M_{\lambda(1)},
\dots, M_{\lambda(m)}$ for which there exists a long exact
sequence $M_{\lambda(1)} \to \cdots \to M_{\lambda(m)}$. When
$m=3$, we recover Fulton's \cite{F2} results on the possible
eigenvalues of majorized Hermitian matrices.
\end{abstract}

\maketitle

\section{Introduction}\label{intro} In \cite{Fr}, Friedland asked for a description of the possible
eigenvalues of Hermitian matrices $A, B$, and $C$ such that $B
\leq A+C$ (i.e., $A+C-B$ is positive semi-definite). A complete
answer to this majorization problem was obtained by Fulton in
\cite{F2} who showed that the eigenvalues of $A, B$, and $C$ are
given by the same inequalities as in Klyachko's theorem \cite{Kl}
for the case when $B=A+C,$ except that the equality
$\tr(B)=\tr(A)+\tr(C)$ is replaced by the linear homogeneous
inequality $\tr(B) \leq \tr(A)+\tr(C).$ As explained in \cite{F2},
the problem about the existence of short exact sequences of finite
abelian $p$-groups \emph{without} zeros at the ends has the exact
same answer as the majorization problem above. In this paper, we
find necessary and sufficient inequalities for the existence of
\emph{long} exact sequences, generalizing Fulton's result.

For every partition $\lambda=(\lambda_1, \dots, \lambda_n)$ and a
(fixed) prime number $p,$ one can construct a finite abelian
$p$-group $M_{\lambda}=\ZZ/p^{\lambda_1}\times\cdots \times
\ZZ/p^{\lambda_n}.$ It is known that every finite abelian
$p$-group is isomorphic to $M_{\lambda}$ for a unique partition
$\lambda.$ Such a group $M_{\lambda}$ is said to be of type
$\lambda.$

For an integer $n \geq 1$, let
$$\mathcal P_n = \{(\lambda_1, \dots, \lambda_n) \in \ZZ^n \mid
\lambda_1 \geq \dots \geq \lambda_n \geq 0 \}
$$
be the semigroup of all partitions with at most $n$ non-zero
parts. Let $m \geq 3$ be a positive integer. We are interested in
the set
$$\Sigma(n, m)=\{ (\lambda(1),
\dots, \lambda(m)) \in \mathcal P_n^m \mid \exists ~M_{\lambda(1)}
\to M_{\lambda(2)} \to \dots \to M_{\lambda(m)} \}.$$ The convex
cone (in $\RR^{nm}$) generated by $\Sigma(n,m)$ is denoted by
$\mathcal C(n,m).$ Now, we are ready to state our first result:
\begin{theorem} \label{satthm.swo} Let $m \geq 3$ and $n \geq 1$ be
two integers.
\begin{enumerate}
\renewcommand{\theenumi}{\arabic{enumi}}
\item The set $\Sigma(n,m)$ is a finitely generated subsemigroup
of $\ZZ^{nm}$ and is saturated, i.e., for every integer $r \geq
1,$
$$
(\lambda(1), \dots, \lambda(m)) \in \Sigma(n,m)
\Longleftrightarrow (r\lambda(1), \dots, r\lambda(m)) \in
\Sigma(n,m).
$$
\item $C(n,m)$ is a rational convex polyhedral cone and
$$
\dim \mathcal C(n,m) = nm.
$$
\end{enumerate}
\end{theorem}

When $m$ is odd, we obtain a recursive method for describing the cone
$\mathcal C(n,m).$ For this, we need to recall some of the terminology from
\cite{CC1}. Let $\lambda(i), 1 \leq i \leq m,$ be $m$ partitions.
Then the \emph{generalized Littlewood-Richardson coefficient}
$f(\lambda(1),\dots, \lambda(m))$ is defined by
$$f(\lambda(1), \dots, \lambda(m))
=\sum c_{\lambda(1),\mu(1)}^{\lambda(2)}\cdot
c_{\mu(1),\mu(2)}^{\lambda(3)}\cdots
c_{\mu(m-4),\mu(m-3)}^{\lambda(m-2)}\cdot
c_{\mu(m-3),\lambda(m)}^{\lambda(m-1)},
$$
where the sum is taken over all partitions $\mu(1), \dots,
\mu(m-3).$ The convention is that when $m=3,$
$f(\lambda(1),\lambda(2),\lambda(3))$ is the Littlewood-Richardson
coefficient $c_{\lambda(1),\lambda(3)}^{\lambda(2)}.$

We refer to the notation paragraph at the end of this section for
the details of our notations. Now, let $(I_1, \dots, I_m)$ be an
$m$-tuple of subsets of $\{1, \dots, n \}$ such that at least one
of them has cardinality at most $n-1.$ We define the following
weakly decreasing sequences of integers (using conjugate
partitions):
$$
\underline\lambda(I_1)=\lambda'(I_1),\hspace{15pt}\underline\lambda(I_m)=
\lambda'(I_m)
$$
and for $2\leq i\leq m-1$
$$\underline\lambda(I_i)=
\begin{cases}
\lambda'(I_i) & \text{if $i$ is even} \\
\lambda'(I_i)-((|I_i|-|I_{i+1}|-|I_{i-1}|)^{n-|I_{i}|}) & \text{if
$i$ is odd}.
\end{cases}
$$

Let $\mathcal S(n,m)$ be the set of all $m$-tuples $(I_1, \dots,
I_m)$ for which:
\begin{enumerate}
\item at least one of the $I_i$ has cardinality at most $n-1;$

\item $|I_1| = |I_2|,$ $|I_{m-1}| = |I_m|;$

\item $\underline \lambda(I_1), \dots, \underline \lambda(I_m)$
are partitions; \item the generalized Littlewood-Richardson
coefficient
$$
f(\underline\lambda(1), \dots, \underline\lambda(m))=1.
$$
\end{enumerate}

For example, if $m=3$ then $\mathcal S(n,3)$ consists of all those
triples $(I_1,I_2,I_3)$ of subsets of $\{1, \dots, n\}$ of the same cardinality $r$ with $r< n$
and
$$
c_{\lambda(I_1), \lambda(I_3)}^{\lambda(I_2)}=1.
$$

The set $\mathcal S(n,m)$ has been used in \cite{CC1} to construct
necessary and sufficient Horn type inequalities for the existence
of long exact sequences of finite abelian $p$-groups with zeros at
the ends. As we are going to see, the same set can be used to
describe $\mathcal C(n,m):$

\begin{theorem} \label{mainthm.swo}
Assume that $m\geq 3$ is odd and let $\lambda(1), \dots,
\lambda(m)$ be $m$ weakly decreasing sequences of $n$ non-negative
real numbers. Then the following are equivalent:
\begin{enumerate}
\renewcommand{\theenumi}{\arabic{enumi}}
\item $(\lambda(1), \dots, \lambda(m)) \in \mathcal C(n,m);$

\item the numbers $\lambda(i)_j$ satisfy
$$
\sum_{i~even}|\lambda(i)|\leq \sum_{i~odd}|\lambda(i)|,$$ and
$$\sum_{i~even} \left(\sum_{j \in
I_i}\lambda(i)_j\right) \leq \sum_{i~odd} \left(\sum_{j\in
I_i}\lambda(i)_j\right)$$ for every $(I_1, \dots, I_m) \in
\mathcal S(n,m);$ if $m > 3$ we also have $$(\lambda(2), \dots,
\lambda(m-1)) \in \mathcal C(n,m-2).$$
\end{enumerate}
\end{theorem}

We should point out that the above theorem fails if $m$ is even
(see Example \ref{counterex}). Nonetheless, for arbitrary $m,$ a
similar description of the cone $C(n,m)$ can be found in Theorem
\ref{mainthm1.swo}.

The strategy for proving the main results of this paper is to show
first that the existence of long exact sequences of finite abelian
$p$-groups without zeros at the ends is equivalent to the
existence of non-zero semi-invariants for a certain quiver. Next,
we use methods from quiver invariant theory developed by Derksen
and Weyman \cite{DW1}, \cite{DW2} to prove Theorem
\ref{satthm.swo} and to find the Horn type inequalities of Theorem
\ref{mainthm.swo} and Theorem \ref{mainthm1.swo}.

The paper is organized as follows. In Section \ref{prelimsec}, we
recall some well-known facts about semi-invariants of quivers and
introduce the cone of effective weights of quivers without
oriented cycles. The quiver setting corresponding to our problem
is defined in Section \ref{quiversettingsec} where we prove
Theorem \ref{satthm.swo}. In Section \ref{hornsec}, we give a
first description of the cone $\mathcal C(n,m)$ and prove Theorem
\ref{mainthm1.swo}. The proof of Theorem \ref{mainthm.swo} is
given in Section \ref{recursec}.

\textbf{Notations.} For a partition $\lambda,$ we denote by
$\lambda'$ the partition conjugate to $\lambda,$ i.e., the Young
diagram of $\lambda'$ is the Young diagram of $\lambda$ reflected
in its main diagonal. We will often refer to partitions as Young
diagrams. If $\lambda=(\lambda_1, \dots, \lambda_N)$ is a weakly
decreasing sequence then we define $r\lambda$ by
$r\lambda=(r\lambda_1, \dots, r\lambda_N).$ Let $\lambda =
(\lambda_1, \dots, \lambda_N)$ and $\mu = (\mu_1, \dots, \mu_M)$
be two sequences of integers. Then we define the sum $\lambda +
\mu$ by first extending $\lambda$ or $\mu$ with zero parts (if
necessary) and then we add them componentwise. If $I=\{z_1< \dots
<z_r\}$ is an $r$-tuple of integers then $\lambda(I)$ is defined
by $\lambda(I)=(z_r-r, \dots, z_1-1).$ For $r \geq 0$ and $a$ two
integers, we denote the $r$-tuple $(a, \dots, a)$ by $(a^r).$ A
composition $\underline{a}$ is just a sequence $\underline{a} =
(a_1, \dots, a_n)$ of non-negative integers. For a weakly
decreasing sequence $\mu$ of $n$ integers, $S^{\mu}(V)$ denotes
the irreducible rational representation of $\GL(V)$ with highest
weight $\mu,$ where $V$ is an $n$-dimensional complex vector
space. Let $\lambda(i)=(\lambda(i)_1, \dots, \lambda(i)_n), 1 \leq
i \leq 3,$ be three weakly decreasing sequences of $n$ integers.
Then we define the Littlewood-Richardson coefficient
$c_{\lambda(1),\lambda(3)}^{\lambda(2)}$ to be the multiplicity of
$S^{\lambda(2)}(\CC^n)$ in $S^{\lambda(1)}(\CC^n) \otimes
S^{\lambda(3)}(\CC^n)$, i.e.
$$
c_{\lambda(1),\lambda(3)}^{\lambda(2)}=\dim_{\CC}
\Hom_{\GL_n(\CC)}(S^{\lambda(2)}(\CC^n), S^{\lambda(1)}(\CC^n)
\otimes S^{\lambda(3)}(\CC^n)).
$$
If $\underline{a} = (a_1, \dots, a_n)$ is a composition and
$\lambda = (\lambda_1, \dots, \lambda_n)$ is a partition with at
most $n$ non-zero parts, we define the Kostka number
$K_{\underline{a}, \lambda}$ to be
$$
K_{\underline{a}, \lambda}=
\dim_{\CC} \Hom_{\GL_n(\CC)}(S^{\lambda}(\CC^n), S^{a_1}(\CC^n)
\otimes \dots \otimes S^{a_n}(\CC^n)).
$$

\subsection*{Acknowledgment} I would like to thank William Fulton
for helpful comments on a preliminary version of this work. I am
grateful to my advisor, Harm Derksen, for many enlightening
discussions on the subject.

\section{Preliminaries}
\label{prelimsec}

\subsection{Generalities}A quiver $Q=(Q_0,Q_1,t,h)$
consists of a finite set of vertices $Q_0$, a finite set of arrows
$Q_1$, and two functions $t,h:Q_1 \to Q_0$ that assign to each
arrow $a$ its tail $ta$ and its head $ha,$ respectively. We write
$ta{\buildrel a\over\longrightarrow}ha$ for each arrow $a \in
Q_1$.

For simplicity, we will be working over the field $\CC$ of complex
numbers.  A representation $V$ of $Q$ over $\CC$ is a family of
finite dimensional $\CC$-vector spaces $\lbrace V(x) \mid x \in
Q_0\rbrace$ together with a family $\{ V(a):V(ta)\rightarrow V(ha)
\mid a \in Q_1 \}$ of $\CC$-linear maps. If $V$ is a
representation of $Q$, we define its dimension vector $\underline
d_V$ by $\underline d_V(x)=\dim_{\CC} V(x)$ for every $x\in Q_0$.
Thus the dimension vectors of representations of $Q$ lie in
$\Gamma=\ZZ^{Q_0}$, the set of all integer-valued functions on
$Q_0$. For every vertex $x,$ the dimension vector of the simple
representation corresponding to $x$ is denoted by $e_x$, i.e.,
$e_x(y)=\delta_{x,y}, \forall y\in Q_0,$ where $\delta_{x,y}$ is
the Kronecker symbol.

Given two representations $V$ and $W$ of $Q$, we define a morphism
$\phi:V \rightarrow W$ to be a collection of linear maps $\lbrace
\phi(x):V(x)\rightarrow W(x)\mid x \in Q_0 \rbrace$ such that
$$\phi(ha)V(a)=W(a)\phi(ta),$$ for
every arrow $a\in Q_1.$ We denote by $\Hom_Q(V,W)$ the
$\CC$-vector space of all morphisms from $V$ to $W$. Let $W$ and
$V$ be two representations of $Q.$ We say that $V$ is a
subrepresentation of $W$ if $V(x)$ is a subspace of $W(x)$ for all
vertices $x \in Q_0$ and $V(a)$ is the restriction of $W(a)$ to
$V(ta)$ for all arrows $a \in Q_1.$ In this way, we obtain the
abelian category $\Rep(Q)$ of all quiver representations of $Q.$ A
dimension vector $\beta$ is said to be a \emph{Schur} root if
there exists a $\beta$-dimensional representation $W$ such that
$\End_{Q}(W) = \CC.$

If $\alpha,\beta$ are two elements of $\Gamma$, we define the
Euler form by
\begin{equation}
\langle\alpha,\beta \rangle = \sum_{x \in Q_0}
\alpha(x)\beta(x)-\sum_{a \in Q_1} \alpha(ta)\beta(ha).
\end{equation}

\subsection{Semi-invariants for quivers}
Let $\beta$ be a dimension vector of $Q$. The representation space
of $\beta$-dimensional representations of $Q$ is defined by
$$\Rep(Q,\beta)=\bigoplus_{a\in Q_1}\Hom(\CC^{\beta(ta)}, \CC^{\beta(ha)}).$$
If $\GL(\beta)=\prod_{x\in Q_0}\GL(\beta(x))$ then $\GL(\beta)$
acts algebraically on $\Rep(Q,\beta)$ by simultaneous conjugation,
i.e., for $g=(g(x))_{x\in Q_0}\in \GL(\beta)$ and $V=(V(a))_{a \in
Q_1} \in \Rep(Q,\beta),$ we define $g \cdot V$ by
$$(g\cdot V)(a)=g(ha)V(a)g(ta)^{-1}\ \text{for every}\ a \in Q_1.$$
Note that $\Rep(Q,\beta)$ is a rational representation of the
linearly reductive group $\GL(\beta)$ and the $\GL(\beta)-$orbits
in $\Rep(Q,\beta)$ are in one-to-one correspondence with the
isomorphism classes of $\beta-$dimensional representations of $Q.$

\textbf{From now on, we will assume that our quivers are without
oriented cycles.} Under this assumption, one can show that there
is only one closed $\GL(\beta)-$orbit in $\Rep(Q,\beta)$ and hence
the invariant ring $\text{I}(Q,\beta)= \CC
[\Rep(Q,\beta)]^{\GL(\beta)}$ is exactly the base field $\CC.$

Now, consider the subgroup $\SL(\beta) \subseteq \GL(\beta)$
defined by
$$
\SL(\beta)=\prod_{x \in Q_0}\SL(\beta(x)).
$$
Although there are only constant $\GL(\beta)-$invariant polynomial
functions on $\Rep(Q,\beta)$, the action of $\SL(\beta)$ on
$\Rep(Q,\beta)$ provides us with a highly non-trivial ring of
semi-invariants.

Let $\SI(Q,\beta)= \CC [\Rep(Q,\beta)]^{\SL(\beta)}$ be the ring
of semi-invariants. As $\SL(\beta)$ is the commutator subgroup of
$\GL(\beta)$ and $\GL(\beta)$ is linearly reductive, we have that
$$\SI(Q,\beta)=\bigoplus_{\sigma
\in X^\star(\GL(\beta))}\SI(Q,\beta)_{\sigma},
$$
where $X^\star(\GL(\beta))$ is the group of rational characters of
$\GL(\beta)$ and $$\SI(Q,\beta)_{\sigma}=\lbrace f \in \CC
[\Rep(Q,\beta)] \mid gf= \sigma(g)f, \forall g \in
\GL(\beta)\rbrace$$ is the space of semi-invariants of weight
$\sigma.$ Note that any $\sigma \in \ZZ^{Q_0}$ defines a
rational character of $\GL(\beta)$ by
$$\{g(x) \mid x \in Q_0\} \in \GL(\beta) \mapsto \prod_{x \in
Q_0}(\det g(x))^{\sigma(x)}.$$ In this way, we can identify
$\Gamma=\ZZ ^{Q_0}$ with the group $X^\star(\GL(\beta))$ of
rational characters of $\GL(\beta),$ assuming that $\beta$ is a
sincere dimension vector (i.e. $\beta(x)>0$ for all vertices $x
\in Q_0$). We also refer to the rational characters of
$\GL(\beta)$ as weights.

If $\alpha \in \ZZ^{Q_0}$, we define the weight $\sigma=\langle
\alpha,\cdot \rangle$ by
$$\sigma(x)=\langle \alpha,e_x \rangle~,~\forall x\in
Q_0.$$ Conversely, it is easy to see that for any weight $\sigma
\in \ZZ^{Q_0}$ there is a unique $\alpha \in \ZZ^{Q_0}$ (not
necessarily a dimension vector) such that $\sigma=\langle
\alpha,\cdot \rangle$. Similarly, one can define $\mu = \langle
\cdot,\alpha \rangle$.

\subsection{Derksen-Weyman saturation}
We write $\beta_1\hookrightarrow \beta$ if every
$\beta$-dimensional representation has a
subrepresentation of dimension vector $\beta_1.$ If $\sigma \in
\mathbb R^{Q_0}$ and $\beta \in \ZZ^{Q_0}$ we define
$\sigma(\beta)$ to be
$$
\sigma(\beta)=\sum_{x \in Q_0}\sigma(x)\beta(x).
$$

The \emph{cone of effective weights} associated to $(Q,\beta)$ is
defined by
$$
C(Q,\beta)= \{ \sigma \in \mathbb R^{Q_0}  \mid \sigma(\beta)=0
\text{~and~} \sigma(\beta_1) \leq 0  \text{~for all~} \beta_1
\hookrightarrow \beta \}.$$

Now, let
$$
\Sigma(Q,\beta)=C(Q,\beta)\bigcap \ZZ^{Q_0}$$ be the semigroup of
lattice points of $C(Q, \beta).$ By construction $C(Q,\beta)$ is a
rational convex polyhedral cone and hence $\Sigma(Q,\beta)$ is
saturated and finitely generated.

In \cite{S2}, Schofield constructed semi-invariants of quivers
with remarkable properties. We should point out that these
Schofield semi-invariants have weights of the form $\langle
\alpha, \cdot \rangle$, with $\alpha$ dimension vectors. A
fundamental result due to Derksen and Weyman \cite{DW1} (see also
\cite{SVB}) states that each weight space of semi-invariants is
spanned by Schofield semi-invariants. An important consequence of
this spanning theorem is the following description of $ \Sigma (Q,
\beta)$ (see \cite{DW1}):

\begin{theorem} [Derksen-Weyman saturation]\label{DW-sat}
Let $Q$ be a quiver and let $\beta$ be a sincere dimension vector.
If $\sigma = \langle \alpha, \cdot \rangle \in \ZZ^{Q_0}$ is a
weight with $\alpha \in \ZZ^{Q_0}$ then the following statements
are equivalent:
\begin{enumerate}
\renewcommand{\theenumi}{\arabic{enumi}}

\item  $\sigma \in \Sigma(Q,\beta)$;

\item $\dim \SI(Q, \beta)_{\sigma} \neq 0;$

\item $\alpha$ must be a dimension vector, $\sigma( \beta) = 0$
and $\alpha \hookrightarrow \alpha+\beta.$
\end{enumerate}
In particular, the dimensions of the weight spaces of
semi-invariants are saturated, i.e., if $n \geq 1$ then
$$
\dim \SI(Q, \beta)_{\sigma} \neq 0 \Longleftrightarrow \dim \SI(Q,
\beta)_{n\sigma} \neq 0.
$$
\end{theorem}

We also have the following reciprocity property:

\begin{lemma}\cite[Corollary 1]{DW1}\label{reciprocity}
Let $\alpha$ and $\beta$ be two dimension vectors. Then
$$\dim\SI(Q,\beta)_{\langle \alpha, \cdot \rangle} = \dim
\SI(Q,\alpha)_{-\langle \cdot, \beta \rangle}.$$
\end{lemma}

Now, we can define $(\alpha \circ \beta)$ by
$$(\alpha \circ \beta)=\dim \SI(Q,\beta)_{\langle \alpha,\cdot
\rangle}=\dim \SI(Q,\alpha)_{-\langle \cdot, \beta \rangle}.$$

In case $\beta$ is a Schur root, we have the following refinement
of Theorem \ref{DW-sat} which is also due to Derksen and Weyman
\cite[Corollary 5.2]{DW2}:

\begin{proposition}\label{descfacets}
Let $Q$ be a quiver with $N$ vertices and let $\beta$ be a Schur
root. Then
\begin{enumerate}
\renewcommand{\theenumi}{\arabic{enumi}}
\item $\dim C(Q,\beta)=N-1.$

\item $\sigma \in C(Q,\beta)$ if and only if $\sigma(\beta)=0$ and
$\sigma(\beta_1)\leq 0$ for every decomposition $\beta=c_1\beta_1
+ c_2\beta_2$ with $\beta_1,\beta_2$ Schur roots,~
$\beta_1\circ\beta_2=1$ and $c_i=1$ whenever $\langle \beta_i,
\beta_i \rangle <0.$
\end{enumerate}
\end{proposition}

Finally, we record a theorem of Schofield on Schur roots which
will be used in the proof of Lemma \ref{lemmafacets.swo}. 

\begin{theorem}\cite[Theorem 6.1]{S1}\label{stabschur}
Let $Q$ be a quiver and let $\beta$ be a dimension vector. Then the
following are equivalent:
\begin{enumerate}
\renewcommand{\theenumi}{\arabic{enumi}}
\item $\beta$ is a Schur root; \item $\sigma_{\beta}(\beta')<0,
\forall~ \beta'\hookrightarrow\beta,~\beta'\neq 0, \beta,$ where
$\sigma_{\beta}=\langle \beta,\cdot\rangle-\langle \cdot,\beta
\rangle.$
\end{enumerate}
\end{theorem}


\section{Long exact sequences from semi-invariants}
\label{quiversettingsec}

In this section, we show that the existence of long exact
sequences of finite abelian $p$-groups without zeros at the ends
is equivalent to the existence of semi-invariants of a certain
quiver. To be more precise, let $(Q, \beta)$ be the following
quiver setting:

\begin{enumerate}
\renewcommand{\theenumi}{\arabic{enumi}}
\item the quiver $Q$ has $m + 1$ central vertices denoted by $0,$
$1=(n,1),$ $2=(n,2),$ $\dots,$ $m=(n,m)$ such that at vertices $1,
2, \dots, m$ we attach $m$ equioriented type $\mathbb A_n$ quivers
(call them flags or arms) $\mathcal F(1), \dots, \mathcal F(m)$
with $\mathcal F(i)$ going in the central vertex $i$ if $i$ is
even and going out from the central vertex $i$ if $i$ is odd;
there are $m-1$ main arrows $a_1, \dots, a_{m-1}$ connecting the
central vertices such that $i+1{\buildrel
a_i\over\longrightarrow}i$ if $i$ is odd and $i{\buildrel
a_i\over\longrightarrow}i+1$ if $i$ is even. Furthermore, there
are $n$ arrows going from vertex $0$ to vertex $1$ and there are
$n$ arrows going from $0$ to $m$ if $m$ is odd; the $n$ arrows go
from $m$ to $0$ if $m$ is even. For example, if $m$ is odd, then
the quiver $Q$ looks like:
$$
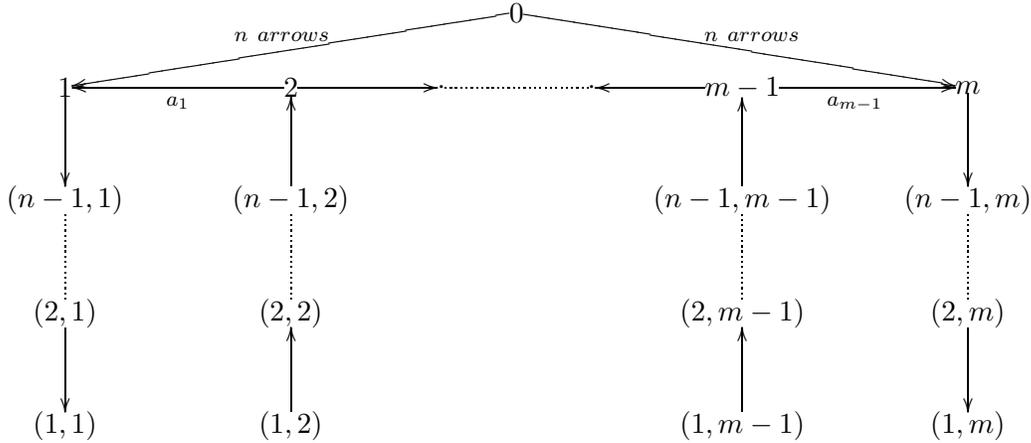
\xy (0, 0)*{0}="a";
        (60, -10)*{m}="a1";
        (60,-25)*{(n-1,m)}="a2";
        (60,-40)*{(2,m)}="a3";
        (60,-55)*{(1,m)}="a4";
        (30,-10)*{m-1}="b1";
        (30,-25)*{(n-1,m-1)}="b2";
        (30,-40)*{(2,m-1)}="b3";
        (30,-55)*{(1,m-1)}="b4";
        (-30,-10)*{2}="c1";
        (-30,-25)*{(n-1,2)}="c2";
        (-30,-40)*{(2,2)}="c3";
        (-30,-55)*{(1,2)}="c4";
        (-60,-10)*{1}="d1";
        (-60,-25)*{(n-1,1)}="d2";
        (-60,-40)*{(2,1)}="d3";
        (-60,-55)*{(1,1)}="d4";
        (-10,-10)*{\cdot}="l";
        (10,-10)*{\cdot}="r";
        {\ar^{n~arrows} "a";"a1"};
        {\ar@{->} "a1";"a2"};
        {\ar@{.} "a2";"a3"};
        {\ar@{->} "a3";"a4"};
        {\ar_{a_{m-1}} "b1";"a1"};
        {\ar@{->} "b2";"b1"};
        {\ar@{.} "b2";"b3"};
        {\ar@{->} "b4";"b3"};
        {\ar^{a_1} "c1";"d1"};
        {\ar@{->} "c2";"c1"};
        {\ar@{.} "c3";"c2"};
        {\ar@{->} "c4";"c3"};
        {\ar_{n~arrows} "a";"d1"};
        {\ar@{->} "d1";"d2"};
        {\ar@{.} "d2";"d3"};
        {\ar@{->} "d3";"d4"};
        {\ar@{->} "c1";"l"};
        {\ar@{->} "b1";"r"};
        {\ar@{.} "l";"r"};
    \endxy
$$



\item the dimension vector $\beta$ is given by $\beta(j, i)=j$ for
all $j \in \{1, \dots, n\},$ $i \in \{1, \dots, m\},$ and
$\beta(0) = 1,$ i.e., $\beta$ is equal to
$$
\begin{matrix}
  &   & 1      &  &\\
n & n & \cdots & n &n\\
n-1 & n-1 & \cdots & n-1 &n-1\\
\vdots& \vdots & &\vdots &\vdots\\
2 & 2 & \cdots & 2 &2\\
1 & 1 & \cdots & 1 &1
\end{matrix}$$
\end{enumerate}

Let $\lambda(1), \dots, \lambda(m)$ be $m$ sequences of $n$ real
numbers. Then we define the weight $\sigma_{\lambda}$ by
\begin{equation} \label{thewt1.swo}
\sigma_{\lambda}(j, i) = (-1)^i(\lambda(i)_j-\lambda(i)_{j+1}),
\forall 1 \leq j \leq n, \forall 1 \leq i \leq m,
\end{equation}
with the convention that $\lambda(i)_{n+1} = 0$ and
\begin{equation} \label{thewt2.swo}
\sigma_{\lambda}(0) = - \sum_{1 \leq i \leq m} \sum_{1 \leq j \leq
n} \sigma_{\lambda}(j, i)j = \sum_{i \text{~odd}}|\lambda(i)| -
\sum_{i \text{~even}}|\lambda(i)|.
\end{equation}

Note that $(\ref{thewt2.swo})$ is equivalent to
$\sigma_{\lambda}(\beta) = 0.$

\begin{lemma}\label{compute.swo}
Let $\lambda(1), \dots, \lambda(m)$ be $m$ partitions with at most
$n$ non-zero parts. Then
$$\dim \SI(Q,\beta)_{\sigma_{\lambda}}\neq 0 \Longleftrightarrow
(\lambda(1), \dots, \lambda(m)) \in \Sigma(n, m).$$
\end{lemma}

\begin{proof} First, we compute the space of semi-invariants
$\SI(Q,\beta)_{\sigma_{\lambda}}$. This is a standard computation
involving Schur functors. For simplicity, let us define $V_j(i)=
\CC^{\beta(j, i)}$. Using the same arguments as in \cite[Lemma
3.1]{CC1}, one can show that each flag $\mathcal F(l)$ going out
from the central vertex $(n,l)$ contributes to $\SI(Q,
\beta)_{\sigma_{\lambda}}$ with
$$
S^{\gamma^{n-1}(l)}V_n(l),
$$
where
$$
\gamma^{n-1}(l) =
((n-1)^{-\sigma_{\lambda}(n-1,l)},\dots,1^{-\sigma_{\lambda}(1,
l)})'.
$$
Now, it is easy to see that
$$
\gamma^{n-1}(l) = (\lambda(l)_1-\lambda(l)_n, \dots,
\lambda(l)_{n-1}-\lambda(l)_n).
$$

Similarly, if $\mathcal F(i)$ is a flag going in the central
vertex $(n,i),$ then its contribution to $\SI(Q,
\beta)_{\sigma_{\lambda}}$ is
$$
S^{\gamma^{n-1}(i)}V^{*}_n(i),
$$
where
$$
\gamma^{n-1}(i) = (\lambda(i)_1-\lambda(i)_n, \dots,
\lambda(i)_{n-1}-\lambda(i)_n).
$$
So far, we have found those spaces of semi-invariants coming from
the vertices of the $m$ flags, except for the central vertices
$i,$ where $i \in \{0,1,\dots, m \}.$ Taking into account the
weights attached to the central vertices,
one can easily see that:
$$
\dim \SI(Q,\beta)_{\sigma_{\lambda}} = \sum
K_{\underline{a},\mu(0)}\cdot c_{\mu(0),\mu(1)}^{\lambda(1)}\cdot
c_{\mu(1),\mu(2)}^{\lambda(2)}\cdot \dots \cdot
c_{\mu(m-1),\mu(m)}^{\lambda(m)}\cdot K_{\underline{b},\mu(m)},
$$
where the sum is over all partitions $\mu(0), \dots,\mu(m)$ and
compositions $\underline{a},\underline{b}$ with
$|\underline{a}|+(-1)^{m+1}|\underline{b}|=|\lambda(1)|-|\lambda(2)|+
\dots+(-1)^{m+1}|\lambda(m)|.$

Now let us prove the implication $"\Rightarrow".$ If $\dim
\SI(Q,\beta)_{\sigma_{\lambda}} \neq 0$
then there exist partitions $\mu(0), \dots,\mu(m)$ such that
$$
f(\mu(0), \lambda(1), \dots, \lambda(m), \mu(m))\neq 0.
$$
This together with Klein's theorem \cite{Kle} imply the existence
of a long exact sequence without zeros at the ends of finite
abelian $p$-groups of types $\lambda(1), \dots,\lambda(m),$ i.e.,
$(\lambda(1), \dots, \lambda(m)) \in \Sigma(n, m).$

For the other implication $"\Leftarrow"$, we extend the given exact
sequence to a long exact sequence with zeros at the ends by taking
the kernel (say, of type $\mu(0)$) of the first morphism and the
cokernel (say, of type $\mu(m)$) of the last morphism of our long
exact sequence. Now, let us break this
long exact sequence with zeros at the ends in short exact
sequences by taking cokernels:
$$0\to M_{\mu(0)}\to M_{\lambda(1)}\to M_{\mu(1)}\to 0,$$ $$0\to M_{\mu(1)}\to M_{\lambda(2)}\to M_{\mu(2)}\to 0,$$
$$\cdots$$ $$0\to M_{\mu(m-1)}\to M_{\lambda(m)}\to M_{\mu(m)}\to 0.$$

Using Klein's theorem \cite{Kle}, this is equivalent to
$$K_{\underline{a},\mu(0)} \cdot c_{\mu(0),\mu(1)}^{\lambda(1)}\cdot
c_{\mu(1),\mu(2)}^{\lambda(2)}\cdots
c_{\mu(m-1),\mu(m)}^{\lambda(m)}\cdot K_{\underline{b},\mu(m)}
\neq 0 ,$$ where $\underline a=\mu(0)$ and $\underline b=\mu(m).$
This implies $\dim\SI(Q,\beta)_{\sigma_{\lambda}} \neq 0.$
\end{proof}

\begin{remark}\label{Rmk-reverse-orient} Note the lemma above remains true when we work with
the quiver obtained from $Q$ by reversing all arrows. Of course,
in this case the new weight is just $-\sigma_{\lambda}$. This
observation is particular useful when proving Lemma
\ref{tech2.swo}.
\end{remark}

\begin{lemma}\label{semigr.swo}
The map
$$
\begin{aligned}
\mathcal C(n,m) & \longrightarrow C(Q, \beta) \\
\lambda = (\lambda(1), \dots, \lambda(m)) & \longrightarrow
\sigma_{\lambda},
\end{aligned}
$$
is an isomorphism of cones that restricts to an isomorphism
between the semigroups of the lattice points.
\end{lemma}

\begin{proof}
The map is well-defined because of Lemma \ref{compute.swo} and the
fact that
$$
\sigma_{\alpha \lambda+\beta \gamma}=\alpha\sigma_{\lambda}+\beta
\sigma_{\gamma},
$$
for all $\alpha, \beta$ (non-negative) real numbers. Note also
that the map is injective. To complete the proof, we only need to
show that the map is surjective.

Let $\sigma \in \Sigma(Q, \beta).$ For $1 \leq j \leq n$ and $1
\leq i \leq m,$ define
$$\beta_1 =
\begin{cases}
\beta - e_{(j, i)} & \text{if $i$ is even} \\
e_{(j, i)} & \text{if $i$ is odd}
\end{cases}
$$
Then it is easy to see that $\beta_1\hookrightarrow \beta$ and
$\sigma(\beta_1)=(-1)^{i+1}\sigma(j, i).$ So, $\sigma$ must
satisfy the so called chamber inequalities, i.e.,
$$
(-1)^i\sigma(j, i) \geq 0,
$$
for all $1 \leq j \leq n$ and $1 \leq i \leq m.$ Now, define $\lambda(i)=(\lambda(i)_1, \dots, \lambda(i)_n)$ by
$$
\lambda(i)_j=(-1)^i\sum_{j \leq k \leq n}\sigma(k, i), \forall 1
\leq i \leq m, 1 \leq j \leq n.
$$
Then the $\lambda(i)$ are partitions with at most $n$ non-zero parts
and $\sigma=\sigma_{\lambda}.$ Hence $(\lambda(1), \dots,
\lambda(m)) \in \Sigma(n,m)$ by Lemma \ref{compute.swo} and this
finishes the proof.

\end{proof}

\begin{lemma}\label{schurmajdim.swo}
The dimension vector $\beta$ is a Schur root of $Q.$
\end{lemma}

\begin{proof}
The dimension vector $\beta$ is in the fundamental region and the
greatest common divisor of its coordinates is one. Then it follows
from \cite[Theorem B(d)]{Kac} that $\beta$ must be a Schur root.
\end{proof}


\begin{proof}[Proof of Theorem \ref{satthm.swo}]
$(1)$ This follows from Derksen-Weyman Saturation Theorem
\ref{DW-sat} and Lemma \ref{semigr.swo}.

$(2)$ As $\beta$ is a Schur root, we know that $\dim C(Q, \beta)$
is the number of vertices of $Q$ minus one and so $ \dim \mathcal
C(n, m) = nm.$
\end{proof}

\section{Horn type inequalities}\label{hornsec}

We work with the quiver set up $(Q, \beta)$ introduced in the
previous section.

\begin{lemma}\label{lemmafacets.swo}
Let $\lambda(1), \dots, \lambda(m)$ be weakly decreasing sequences
of $n$ real numbers. Then
$$
\sigma_{\lambda} \in C(Q, \beta) \Longleftrightarrow
\sigma_{\lambda}(\beta_1) \leq 0,
$$
for every dimension vector $\beta_1 \neq \beta$ with $\beta_1
\circ (\beta - \beta_1) = 1$ and $\beta_1$ weakly increasing with
jumps of at most one along the $m$ flags (from bottom to top).
\end{lemma}

\begin{proof} The implication
$"\Longrightarrow"$ follows from Theorem \ref{DW-sat}{(3)}.

Now, let us prove $"\Longleftarrow".$ Using Proposition
\ref{descfacets} and $\sigma_{\lambda}(\beta)=0$, we only need to show
$$
\sigma_{\lambda}(\beta_1) \leq 0,
$$
for every decomposition $\beta = c_1\beta_1 + c_2\beta_2$ with
$\beta_1,~ \beta_2$ Schur roots and $\beta_1 \circ \beta_2 =1$.

If $\beta_1$ is either $\beta_1 = e_{(j, i)}$ for some $1 \leq j
\leq n-1$ and $i$ odd or $\beta_1 = \beta -e_{(j, i)}$ for some $1
\leq j \leq n-1$ and $i$ even then $\sigma_{\lambda}(\beta_1) \leq
0$ is equivalent to $\lambda(1), \dots, \lambda(m)$ being weakly
decreasing sequences.

Now, let us assume $\beta_1$ is not of the above form. We are
going to show that $c_1 = c_2 = 1$ and that $\beta_1, \beta_2$ are
weakly decreasing with jumps of at most one along the $m$ flags
(from bottom to top). Let us denote $c_1\beta_1=\beta',
c_2\beta_2=\beta''$. Since $\beta'\circ\beta'' \neq 0$ it follows
from Theorem \ref{DW-sat} that any representation of dimension
vector $\beta$ has a subrepresentation of dimension vector
$\beta'.$ Therefore, $\beta'$ must be weakly increasing along each
flag going in and it has jumps of at most one along each flag
going out.

Next, we will show that $\beta'$ has jumps of at most one along
each flag $\mathcal F(i)$ going in a central vertex and $\beta'$
is weakly increasing along each flag $\mathcal F(i)$ going out of
a central vertex. For simplicity, let us write
\begin{diagram} \mathcal F(i):
1&\rTo&2&\cdots&n-1&\rTo&n,
\end{diagram}
for a flag going in its central vertex $(n,i)$ (i.e., $i$ is even).
Assume to the contrary that there is an $l\in \{1, \dots, n-1\}$
such that $\beta'(l+1)>\beta'(l)+1$. Then $\beta ''(l+1) < \beta
''(l)$ which implies that $e_l\hookrightarrow\beta ''$.
Since $\beta''$ is $\langle\beta',\cdot \rangle$-semi-stable it
follows that $\langle \beta',e_l\rangle\leq 0$. So,
$\beta'(l)\leq \beta'(l-1)$ and hence $\beta'(l)=\beta'(l-1)$ or
$\beta''(l)=\beta''(l-1)+1$. This shows that $c_2=1$ and $\beta''
- e_l\hookrightarrow \beta''$. From the fact that $\beta
''(=\beta_2)$ is a Schur root and Theorem \ref{stabschur} we
obtain that $\beta ''$ is $\sigma_{\beta ''}$-stable. Since
$e_l\hookrightarrow\beta
''$, $\beta''-e_l\hookrightarrow\beta ''$ and $\beta ''
\neq e_l$ it follows $\langle\beta
'',e_l\rangle-\langle e_l, \beta''\rangle<0$ and
$\langle \beta '',\beta ''-e_l\rangle-\langle \beta '' -
e_l, \beta '' \rangle < 0$. But this is a contradiction.
We have just proved that $\beta'$ has jumps of at most one along
each flag going in. Similarly, one can show that $\beta'$ has to
be weakly increasing along each flag going out.

Now, let us show that $c_1=c_2=1$. Since $\beta'=c_1\beta_1$ has
jumps of at most one along each flag, we obtain $0\leq
c_1(\beta_1(l+1, i)-\beta_1(l, i))\leq 1$ for all $l\in \{1,
\dots, n-1\}$ and $i \in \{1, \dots, m\}.$ If there are $l,i$ such
that $\beta_1(l+1, i)-\beta_1(l, i)\neq 0$ then $c_1=1$.
Otherwise, there is an $i$ such that $\beta'(1, i)=1$ and so
$c_1=1$. Similarly, one can show $c_2=1$.

In conclusion, $\beta=\beta_1+\beta_2$ with $\beta_1$ weakly
increasing with jumps of at most one along the $m$ flags. So, we have $\sigma_{\lambda}(\beta_1) \leq
0$ and we are done.
\end{proof}

\begin{remark}
We want to point out that some of the inequalities obtained in
Lemma \ref{lemmafacets.swo} are redundant. The reason for the
redundancy is that some of the $\beta_1$ or
$\beta_2=\beta-\beta_1$ above might not be Schur roots.
\end{remark}

\begin{example}\label{casen=1.swo}
Let $n = 1$ and $m \geq 3.$ Let $\lambda(i)= (\lambda_i)$ with
$\lambda_i$ non-negative integers, $1 \leq i \leq m$. We show that
there exists an exact sequence
$$
 \ZZ/ p^{\lambda_1}\to\ZZ/ p^{\lambda_2} \to \dots \to \ZZ/ p^{\lambda_m}
$$
if and only if
$$
\lambda_i-\lambda_{i+1}+\dots-\lambda_{j-1}+\lambda_j \geq 0,
$$
for all even numbers $i$ and $j$ with $ 2 \leq i \leq j \leq m$
and
$$
\lambda_{i'}-\lambda_{i'+1}+\dots-\lambda_{j'-1}+\lambda_{j'} \geq
0,
$$
for all odd numbers $i'$ and $j'$ with $1 \leq i' \leq j' \leq m.$

The quiver we work with in this case is of type
$\tilde{\mathbb{A}}_m.$ For example, if $m$ is odd then the quiver
looks like:

$$
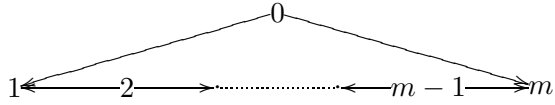
\xy (0, 0)*{0}="a";
        (35, -10)*{m}="a1";
        (20,-10)*{m-1}="b1";
        (-20,-10)*{2}="c1";
        (-35,-10)*{1}="d1";
        (-8,-10)*{\cdot}="l";
        (8,-10)*{\cdot}="r";
        {\ar@{->} "a";"a1"};
        {\ar@{->} "b1";"a1"};
        {\ar@{->} "c1";"d1"};
        {\ar@{->} "a";"d1"};
        {\ar@{->} "c1";"l"};
        {\ar@{->} "b1";"r"};
        {\ar@{.} "l";"r"};
    \endxy
$$
The dimension vector $\beta$ is
$$
\begin{matrix}
  &   & 1      &  &\\
n & n & \cdots & n &n.
\end{matrix}$$

We want to find all Schur roots $\beta_1$ and $\beta_2$ such that
$\beta_1 \neq \beta$ and $\beta_1\circ \beta_2 = 1.$

\emph{Case 1.} If $\beta_1(0) = 1$ then $\beta_1$ has to be of the
form $$\beta_1(v) =
\begin{cases}
0 & \text{if $i \leq v \leq j,$} \\
1 & \text{otherwise},
\end{cases}
$$
for two even numbers $i$ and $j,$ $2 \leq i \leq j \leq m.$
Conversely, any dimension vector $\beta_1$ of this form has the
property that $\beta_1,~ \beta - \beta_1$ are Schur roots and
$\beta_1 \circ (\beta - \beta_1) = 1.$ In this case, we have
$$
\sigma_{\lambda}(\beta_1)=\sum_{\buildrel {i \leq v \leq j} \over
{v~ \text{odd}}} \lambda_v - \sum_{\buildrel {i \leq v \leq
j}\over {v~ \text{even}}}\lambda_v.
$$

\emph{Case 2.} If $\beta_1(0) = 0$ then $\beta_1$ has to be of the
form $$\beta_1(v) =
\begin{cases}
1 & \text{if $i' \leq v \leq j',$} \\
0 & \text{otherwise},
\end{cases}
$$ for two odd numbers $i'$ and $j',$ $1 \leq i' \leq j' \leq m.$ Again, if
$\beta_1$ is of this form then $\beta_1,~ \beta - \beta_1$ are
Schur roots and $\beta_1 \circ (\beta - \beta_1) = 1.$ In this
case, we have
$$
\sigma_{\lambda}(\beta_1)=\sum_{\buildrel {i' \leq v \leq j'}\over
{v~ \text{even}}}\lambda_v - \sum_{\buildrel {i' \leq v \leq
j'}\over {v~ \text{odd}}}\lambda_v.
$$
\end{example}

\vspace{20pt}

In what follows, we find a closed form of  those inequalities
obtained in Lemma \ref{lemmafacets.swo}. Let $\beta_1$ be a
dimension vector which is weakly increasing with jumps of at most
one along the $m$ flags of $Q.$ Define the sets
$$
I_i = \{ l \mid \beta_1(l, i) > \beta_1(l-1, i)),1 \leq l\leq n \}
$$
with the convention that $\beta_1(0, i)=0$ for all $1 \leq i \leq
m.$ Then it is easy to see that $|I_i| = \beta_1(i), \forall 1
\leq i \leq m.$

Conversely, given an $m$-tuple $I = (I_1, \dots, I_m)$ of subsets
of $\{1, \dots, n \},$ we can construct two dimension vectors
$\beta_{I}$ and $\beta'_I$ as follows. If
$$
I_i = \{z(i)_1< \dots < z(i)_r\},
$$
we define
$$
\beta_I(k, i) = \beta'_I (k, i) = j-1, \forall z(i)_{j-1}\leq k <
z(i)_j, \forall 1 \leq j \leq r+1,
$$
with the convention that $z(i)_0=0$ and $z(i)_{r+1}=n+1$ for all
$1 \leq i \leq m.$ At vertex $0,$ we let $\beta_{I}(0) = 0$ and
$\beta'_I(0) = 1.$

\begin{theorem}\label{mainthm1.swo}
The cone $\mathcal C(n,m)$ consists of all $m$-tuples
$(\lambda(1), \dots, \lambda(m))$ of weakly decreasing sequences
of $n$ real numbers for which:
\begin{enumerate}
\renewcommand{\theenumi}{\arabic{enumi}}
\item

$$\sum_{i \text{~even}}\left( \sum_{j \in
I_i}\lambda(i)_j\right) \leq \sum_{i \text{~odd}}\left( \sum_{j
\in I_i}\lambda(i)_j\right),$$ for every $m$-tuple $I = (I_1,
\dots, I_m)$ of subsets of $\{1, \dots, n \}$ with $$\beta_I \circ
(\beta - \beta_I) = 1;$$ \item

$$\sum_{i \text{~odd}} \left( \sum_{j \notin I_i}\lambda(i)_j\right) \leq \sum_{i
\text{~even}} \left( \sum_{j \notin I_i}\lambda(i)_j \right),$$
for every $m$-tuple $I = (I_1, \dots, I_m)$ of subsets of $\{1,
\dots, n \}$ with $$\beta'_I \circ (\beta - \beta'_I) = 1.$$
\end{enumerate}
\end{theorem}

\begin{proof}
From Lemma \ref{semigr.swo} it follows that
$$
(\lambda(1), \dots, \lambda(m)) \in \mathcal{C}(n,m)
\Longleftrightarrow \sigma_{\lambda} \in C(Q, \beta).
$$

Now, let $\beta_1$ be a dimension vector which is weakly
increasing with jumps of at most one along the $m$ flags, $\beta_1
\neq \beta$  and $\beta_1 \circ (\beta - \beta_1) = 1.$ Let $I =
(I_1, \dots, I_m)$ be the jump sets. Then $\beta_1$ is $ \beta_I$
if $\beta_1(0)=0$ or $\beta_1$ is $\beta'_I$ if $\beta_1(0)=1.$
Moreover, we have that
$$
\sigma_{\lambda}(\beta_I) = \sum_{i~even} \left(\sum_{j \in
I_i}\lambda(i)_j\right) - \sum_{i~odd} \left(\sum_{j\in
I_i}\lambda(i)_j\right),
$$
$$
\sigma_{\lambda}(\beta'_I) = \sum_{i~odd} \left(\sum_{j \notin
I_i}\lambda(i)_j\right)-\sum_{i~even} \left(\sum_{j\notin
I_i}\lambda(i)_j\right)
$$
and, of course, $\sigma_{\lambda}(\beta) = 0.$ The proof follows
now from Lemma \ref{lemmafacets.swo}.
\end{proof}

\begin{remark} It is easy to see that if $\lambda(i)$ are weakly
decreasing sequences satisfying the conditions $(1)$ and $(2)$ of
Theorem \ref{mainthm1.swo} then $\lambda(i)$ are sequences of
non-negative real numbers. Of course, this non-negativity is
automatically satisfied in $C(n,m).$
\end{remark}

\section{A recursive description}
\label{recursec}

First, we recall a reduction method that appears in \cite{DW1},
\cite{DL}, \cite{SW2}, and \cite{SW1}.

\begin{lemma}\label{reductionlemma.swo}
Let $Q$ be a quiver and $v_0$ a vertex such that near $v_0,$ $Q$
looks like:
$$
v_1 \rTo{a} v_0 \rTo{b} w_1.
$$
Suppose that $\beta$ is a dimension vector and $\sigma$ is a
weight such that $$\beta(v_0)\geq \min \{ \beta(w_1), \beta(v_1)\}
\text{~and~}\sigma(v_0)=0.$$ Let $\overline Q$ be the quiver
defined by $\overline Q_0 = Q_0 \setminus \{ v_0 \}$ and
$\overline Q_1 = (Q_1\setminus \{ a, b \}) \cup \{ba \}.$ If
$\overline \beta=\beta|_{\overline Q}$ is the restriction of
$\beta$ and $\overline\sigma=\sigma|_{\overline Q}$ is the
restriction of $\sigma$ to $\overline Q$ then
$$ \SI(Q,\beta)_{\sigma}\cong \SI(\overline
Q,\overline\beta)_{\overline\sigma}.
$$
\end{lemma}


\textbf{From now on we will assume that $m$ is odd.} Under this
assumption, we are able to further describe $\beta_I \circ
(\beta-\beta_I)$ and $\beta'_I\circ (\beta-\beta'_I).$ For the
convenience of the reader, we recall some of the notations from
Section \ref{intro}. Let $(I_1, \dots, I_m)$ be an $m$-tuple of
subsets of $\{1, \dots, n \}$ such that at least one of them has
cardinality at most $n-1.$ We define the following weakly
decreasing sequences of integers (using conjugate partitions):
$$
\underline\lambda(I_1)=\lambda'(I_1),~\underline\lambda(I_m)=
\lambda'(I_m)
$$
and for $2\leq i\leq m-1$
$$\underline\lambda(I_i)=
\begin{cases}
\lambda'(I_i) & \text{if $i$ is even} \\
\lambda'(I_i)-((|I_i|-|I_{i+1}|-|I_{i-1}|)^{n-|I_{i}|}) & \text{if
$i$ is odd}
\end{cases}
$$

\begin{lemma} \label{tech1.swo}
Let $I = (I_1, \dots, I_m)$ be an $m$-tuple of subsets of $\{1,
\dots, n\}$ as above and such that $|I_1| = |I_2|$ and $|I_{m-1}|
= |I_m|.$ If $\beta_I \circ (\beta - \beta_I) \neq 0$ then
$\underline\lambda(I_i)$ are partitions and
$$\beta_I \circ (\beta - \beta_I) = f(\underline\lambda(I_1), \dots,
\underline\lambda(I_m)).$$

Consequently,
$$
\beta_I\circ (\beta-\beta_I)=1 \Longleftrightarrow I \in \mathcal
S(n,m).
$$
\end{lemma}

\begin{proof}
Let us denote $\beta_I$ by $\beta_1$ and $\beta - \beta_I$ by
$\beta_2.$ Then we have that
$$\beta_1 \circ \beta_2 = \dim \SI(Q, \beta_1)_{-\langle \cdot,
\beta_2 \rangle}.$$ Since $\beta_1(0) = 0,$ we can work with the
quiver $Q'$ obtained from $Q$ by deleting the vertex $0$ and all
the arrows going out from this vertex. If $\beta'_1$ and
$\beta'_2$ are the restrictions of $\beta_1$ and $\beta_2$ to
$Q',$ then the restriction of the weight $-\langle \cdot, \beta_2
\rangle$ to $Q'$ is exactly $-\langle \cdot, \beta'_2 \rangle$ as
the $n$ arrows connecting vertex $0$ and $m$ point towards vertex
$m.$ Therefore, we have
$$\beta_1 \circ \beta_2 = \beta'_1 \circ \beta'_2.$$ Let us denote $
\langle \beta'_1, \cdot \rangle$ by $\sigma'_1$. As $\beta'_1(1) =
\beta'_1(2) = |I_1| = |I_2|$ and $\beta'_1(m-1) = \beta'_1(m) =
|I_{m - 1}| = |I_m|$ it follows that $\sigma'_1(1)= \sigma'_1(m) =
0.$

At this point, we can apply the reduction Lemma
\ref{reductionlemma.swo} to reduce $Q'$ to the quiver $Q''$
obtained from $Q'$ by removing the two vertices $1$ and $m.$
Again, it easy to check that if $\beta''_1,~\beta''_2$ are the
restriction of $\beta'_1,~\beta'_2$ to $Q''$ then $$\beta'_1 \circ
\beta'_2 = \beta''_1 \circ \beta''_2.$$ On the other hand, this
reduced quiver $Q''$ is exactly the generalized flag quiver from
\cite[Section 3]{CC1}. It follows from (\cite[Lemma 6.4]{CC1})
that $\underline{\lambda}(i),~ 1 \leq i \leq m$ are partitions and
$$\beta''_1 \circ \beta''_2=f(\underline\lambda(I_1), \dots,
\underline\lambda(I_m)).$$ This finishes the proof.
\end{proof}

\begin{remark} \label{mainrmk.swo}
Let $\beta=\beta_1+\beta_2$ with $\beta_1$ weakly increasing with
jumps of at most one along the flags and $\beta_1\circ\beta_2\neq
0.$ We claim that
$$\beta_1(0)=1 \Rightarrow \beta_1 \text{~is~}\beta \text{~
along the flags~} \mathcal F(1) \text{~and~} \mathcal F(m).$$
Indeed, we have that $\beta_1\hookrightarrow\beta$ by Theorem
\ref{DW-sat}{(3)}. Consider a representation $W \in \Rep(Q,\beta)$
with $\{\Ima W(a_i)\}_{1\leq i\leq n}$ linearly independent. Since
$W$ must have a $\beta_1$-dimensional subrepresentation, we obtain
that $\beta_1(1)=n$ and so $\beta_1$ has to be $\beta$ along
$\mathcal F(1).$ Similarly, as $m$ is odd, we have that
$\beta_1(0) = 1$ implies that $\beta_1$ equals $\beta$ along the
flag $\mathcal F(m).$
\end{remark}

\begin{lemma} \label{tech2.swo}
Let $I = (I_1, \dots, I_m)$ be an $m$-tuple of subsets of $\{1,
\dots, n \}$ and let $\lambda(i), 1 \leq i \leq m$ be weakly
decreasing sequences of $n$ non-negative reals.
\begin{enumerate}
\renewcommand{\theenumi}{\arabic{enumi}}
\item If $ \beta'_I \circ (\beta - \beta'_I) \neq 0 $ and
$(\lambda(2), \dots, \lambda(m-1)) \in \mathcal C(n, m-2)$ then
$$
\sigma_{\lambda}(\beta'_I) \leq 0.
$$
\item Suppose that at least one of the sets $I_1, \dots, I_m$ has
cardinality at most $n-1$ and $\beta_I \circ (\beta - \beta_I) =
1.$ Furthermore, assume that
$$\sum_{i~even} \left(\sum_{j \in J_i}\lambda(i)_j\right) \leq
\sum_{i~odd} \left(\sum_{j\in J_i}\lambda(i)_j\right),$$ for every
$(J_1, \dots, J_m) \in \mathcal S(n,m).$ Then
$$
\sigma_{\lambda}(\beta_I) \leq 0.
$$
\end{enumerate}
\end{lemma}

\begin{proof}
$(1)$ Let us write $\beta_1 = \beta'_I$ and $\beta_2 = \beta -
\beta'_I.$ As $\beta_1(0) = 1$ and $m$ is odd it follows from
Remark \ref{mainrmk.swo} that $\beta_1$ has to be equal to $\beta$
along the flags $\mathcal F(1)$ and $\mathcal F(m).$ In other
words,  $\beta_2$ is zero at vertex $0$ and at all vertices of
the flags $\mathcal F(1)$ and $\mathcal F(m).$

Now, let $Q'$ be the quiver obtained from $Q$ by deleting the
vertex $0$, the flags $\mathcal F(1)$ and $\mathcal F(m)$ and all
the arrows connected with these deleted vertices. If $\beta'_i$ is
the restriction of $\beta_i$ to $Q',$ $i \in \{1, 2 \}$, then
$$\beta_1 \circ \beta_2 = \beta'_1 \circ \beta'_2.$$

Let $Q''$ be the quiver obtained from $Q'$ by adding a new vertex
$0,$ $n$ arrows from vertex $2$ to $0$ and $n$ arrows from vertex
$m - 1$ to $0.$ We denote by $\beta''_1$ and $\beta''_2$ the
extensions of $\beta'_1$ and $\beta'_2$ to $Q''$ such that
$\beta''_1(0) = 1$ and $\beta''_2(0) = 0.$ Again, it easy to see
that $$\beta''_1 \circ \beta''_2 = \beta'_1 \circ \beta'_2.$$ Note
that $Q''$ is the quiver corresponding to $\Sigma(n, m-2),$ except
that all the arrows have the opposite orientation. So, let us
define the weight $\sigma''_{\lambda}$ for $Q''$ by
$$
\sigma''_{\lambda}(j, i) = (-1)^i(\lambda(i)_j-\lambda(i)_{j+1}),
\forall 1 \leq j \leq n, \forall 2 \leq i \leq m-1,
$$
and $\sigma''_{\lambda}(0)$ is determined by
$\sigma''_{\lambda}(\beta'') = 0,$ where $\beta''$ is just the
restriction of $\beta$ to $Q''_0.$

From Remark \ref{Rmk-reverse-orient}, we deduce that
$\sigma''_{\lambda} \in C(Q'', \beta'')$ if and only if
$(\lambda(2), \dots, \lambda(m-1)) \in \mathcal C(n, m-2).$ As
$\beta''_1 \hookrightarrow \beta''$ and $\sigma''_{\lambda} \in
C(Q'', \beta'')$ it follows that $\sigma''_{\lambda}(\beta''_1)
\leq 0,$ i.e.,
$$\sum_{\buildrel {2 \leq i \leq m - 1}\over { i~\text{even}}} \left( \sum_{j \in
I_i}\lambda(i)_j\right) - \sum_{\buildrel {2 \leq i \leq m - 1}
\over {i~\text{odd}}} \left( \sum_{j \in I_i}\lambda(i)_j\right) +
\sum_{\buildrel {2 \leq i \leq m-1} \over
{i~\text{odd}}}|\lambda(i)| - \sum_{\buildrel {2 \leq i \leq m -
1}\over {i~\text{even}}}|\lambda(i)| \leq 0.$$ In other words, we
have
$$\sigma_{\lambda}(\beta'_I) \leq 0.$$

$(2)$ Let $\alpha_1 = \beta_I$ and $\alpha_2 = \beta - \beta_I.$
Again, as $\alpha_1(0) = 0,$ we can simplify our quiver by
deleting the vertex $0$ and all the arrows going out from this
vertex. We denote the simplified quiver by $\widetilde{Q}$ and the
restriction of the dimension vectors will be noted by
$\widetilde{\alpha}_1$, $\widetilde{\alpha}_2$, and
$\widetilde{\beta}$.

Next, we compute the dimension $$\beta_I \circ
(\beta-\beta_I)=\widetilde{\alpha}_1 \circ \widetilde{\alpha}_2 =
\dim \SI(\widetilde{Q}, \widetilde{\alpha}_2)_{\langle
\widetilde{\alpha}_1, \cdot \rangle}$$ using the same arguments as
in Lemma \ref{tech1.swo}. Note that the weight
$\widetilde{\sigma}_1 = {\langle \widetilde{\alpha}_1, \cdot
\rangle}$ is equal to $\widetilde{\alpha}_1(1) -
\widetilde{\alpha}_1(2)$ at vertex $1$ and it is equal to
$\widetilde{\alpha}_1(m) - \widetilde{\alpha}_1(m - 1)$ at vertex
$m.$ Furthermore, as $\widetilde{\alpha}_1 \circ
\widetilde{\alpha}_2 \neq 0,$ we have $\widetilde{\alpha}_1(1)\geq
\widetilde{\alpha}_1(2)$ and $\widetilde{\alpha}_1(m)\geq
\widetilde{\alpha}_1(m-1).$ To see this, just take
$\widetilde{W}\in \Rep(\widetilde{Q}, \widetilde{\beta})$ to be
bijective along the main arrows $a_1$ and $a_{m-1}.$

Note that $I_1, \dots, I_m$ are the jump sets of
$\widetilde{\alpha}_1$ along the $m$ flags of $Q.$ Let $J_1$ be
the subset of $I_1$ consisting of the first
$\widetilde{\alpha}_1(2)$ elements of $I_1.$ Similarly, let $J_m$
be the subset of $I_m$ consisting of the first
$\widetilde{\alpha}_1(m-1)$ elements of $I_m$. As
$\widetilde{\alpha}_1 \circ \widetilde{\alpha}_2 \neq 0$, we know
that $\underline \lambda(J_1),~\underline \lambda(J_m),$
$\underline\lambda(I_i)$ must be partitions for all $2 \leq i \leq
m-1$ and
$$\widetilde{\alpha}_1 \circ \widetilde{\alpha}_2=f(\underline \lambda(J_1),
\underline\lambda(I_2),\dots,\underline\lambda(I_{m-1}),\underline\lambda(J_m)).
$$
It is clear that at least one of the $J_1, I_2, \dots, I_{m-1},
J_m$ has cardinality at most $n-1,$ and hence, $(J_1, I_2, \dots,
I_{m-1}, J_m) \in S(n, m).$ Therefore, we have
$$
\sum_{i~even} \left(\sum_{j \in I_i}\lambda(i)_j\right) \leq
 \sum_{j\in
J_1}\lambda(1)_j + \sum_{j\in J_m}\lambda(m)_j + \sum_{\buildrel
{2 \leq i \leq m-1} \over {i~\text{odd}}} \left(\sum_{j\in
I_i}\lambda(i)_j\right).$$ As $\lambda(1)_j$ and $\lambda(m)_j$
are assumed to be non-negative for all $1 \leq j \leq n$ we obtain
that $\sigma_{\lambda}(\beta_I) \leq 0.$
\end{proof}


\begin{proof}[Proof of Theorem \ref{mainthm.swo}]
First, let us prove that $(1) \Rightarrow (2).$ If $I = (I_1,
\dots, I_m)$ is an $m$-tuple in $\mathcal S(n,m)$ then $\beta_I
\circ (\beta - \beta_I) \neq 0,$ by Lemma \ref{tech1.swo} and so
$\beta_I \hookrightarrow \beta.$ As $\sigma_{\lambda} \in C(Q,
\beta),$ we have that $\sigma_{\lambda}(\beta_I)\leq 0$ which is
equivalent to
$$
\sum_{i~even} \left(\sum_{j \in I_i}\lambda(i)_j\right) \leq
\sum_{i~odd} \left(\sum_{j\in I_i}\lambda(i)_j\right).
$$

To obtain the first inequality, we just note that $\beta - e_{0}
\hookrightarrow\beta$ (this is not true if $m$ is even) and this
clearly implies that
$$
\sum_{i~even}|\lambda(i)|\leq \sum_{i~odd}|\lambda(i)|.
$$
Next, it is clear that $(\lambda(1), \dots, \lambda(m)) \in
\mathcal C(n,m)$ implies $(\lambda(2), \dots, \lambda(m-1)) \in
\mathcal C(n,m-2).$

For the other implication $(1)\Leftarrow (2),$ let $I = (I_1,
\dots, I_m)$ be an $m$-tuple of subsets of $\{1, \dots, n \}.$ If
$|I_i| = n, \forall 1 \leq i \leq m$ then $\beta_I = \beta - e_0$
and $\beta'_I = \beta.$ In this case, we have
$$
\sigma_{\lambda}(\beta_I) = \sum_{i~even}|\lambda(i)| -
\sum_{i~odd}|\lambda(i)| \leq 0,
$$
and $\sigma_{\lambda}(\beta'_I) = 0.$

Now, let us assume that at least one of the $I_i$ has cardinality
at most $n-1.$ If $\beta'_I \circ (\beta - \beta'_I) = 1$ then
$\sigma_{\lambda}(\beta'_I) \leq 0$ by Lemma \ref{tech2.swo}{(1)}.
If $\beta_I \circ (\beta - \beta_I) = 1$ then it follows from
Lemma \ref{tech2.swo}{(2)} that $\sigma_{\lambda}(\beta_I) \leq
0.$ The proof follows now from Theorem \ref{mainthm1.swo}.
\end{proof}





\begin{remark} \label{counterex}
Let us point out that Theorem \ref{mainthm.swo} fails if $m$ is
even. For example, one can take $m = 4,$ $n = 1.$ Then $\lambda(1)
= (3),$ $\lambda(2) = (3),$ $\lambda(3) = (1),$ $\lambda(4) = (2)$
give a counterexample
to Theorem \ref{mainthm.swo}. 
\end{remark}

When $m=3$ in Theorem \ref{mainthm.swo}, we recover Fulton's
result \cite{F2}:

\begin{corollary}[\textbf{Majorization problem}]
Let $\lambda(1), \lambda(2), \lambda(3)$ be three partitions with
at most $n$ non-zero parts. Then the following are equivalent:
\begin{enumerate}
\renewcommand{\theenumi}{\arabic{enumi}}
\item there exist a short exact sequence of the form
$$M_{\lambda(1)} \to M_{\lambda(2)} \to M_{\lambda(3)},$$
where $M_{\lambda(i)}$ is a finite abelian $p$-group of type
$\lambda(i);$

\item the numbers $\lambda(i)_j$ satisfy $$|\lambda(2)| \leq |\lambda(1)|+|\lambda(3)|$$ and
$$
\sum_{j \in I_2}\lambda(2)_j \leq \sum_{j \in
I_1}\lambda(1)_j+\sum_{j \in I_3}\lambda(3)_j
$$
for all triples $(I_1, I_2, I_3)$ of subsets of $\{1, \dots, n\}$
of the same cardinality $r$ with $r<n$ and
$c_{\lambda(I_1),\lambda(I_3)}^{\lambda(I_2)}=1.$
\end{enumerate}
\end{corollary}

\end{document}